\documentclass[12pt]{article}

\usepackage{amsfonts}
\usepackage{amsmath}

\newtheorem{theorem}{Theorem}
\newtheorem{lemma}{Lemma}
\newtheorem{remark}{Remark}
\newenvironment{proof}{\noindent{\it Proof. }\rm}
{\unskip\nobreak\hfil\penalty50\hskip1em\hbox{}
\nobreak\hfill\qed\par\smallskip}
\def\qed{\vrule height1ex width1ex depth0pt}

\begin{document}
\title{On representations of a q-analogue of the $*$-algebra $Pol(Mat_{2,2})$}
\author{Lyudmila  Turowska \vspace{0.3cm}\\
{\footnotesize \sl Department of Mathematics, Chalmers University of 
Technology, }\\ 
{\footnotesize \sl 412 96 G\"oteborg, Sweden} }

\date{}
\maketitle{}
\abstract{Bounded Hilbert space $*$-representations are studied for a 
$q$-analogue of the $*$-algebra $Pol(Mat_{2,2})$ of polynomials on the space 
$Mat_{2,2}$ of complex $2\times 2$ matrices.}

\bigskip

\noindent
{\bf 1. Introduction}
\bigskip

The study of q-analogues of the Cartan domains (irreducible bounded symmetric 
domains) was initiated by S.~Sinel'shchikov and L.~Vaksman  in \cite{SV}. 
In particular, for each Cartan domain they defined the $*$-algebra 
$Pol({\mathfrak g}_{-1})_q$, a 
q-analogue of the  polynomial algebra on the prehomogenous vector space 
${\mathfrak g}_{-1}$, and set a problem on
 investigation of their representations. The theory of representations of the 
$*$-algebras corresponding to domains of rank 1 is well-understood. In this 
paper  our purpose is to study such representations for one of the 
popular Cartan domains of rank 2, the matrix ball in the space
$Mat_{2,2}$ of complex $2\times2$  matrices. Following
\cite{SSV} we will denote this $*$-algebra by $Pol(Mat_{2,2})_q$. 
A description of $Pol(Mat_{m,n})_q$, $m$, $n\in {\mathbb N}$, in  
 terms of generators and relations is given in \cite{SSV}. In the paper
we classify all irreducible representations of $Pol(Mat_{2,2})_q$ by bounded
operators on a  Hilbert space. The method which we use here is based
on the study of some dynamical system arising on a spectrum of a commutative
$*$-subalgebra of $Pol(Mat_{2,2})_q$ (see 
\cite{OS}). Note that the $*$-algebra has also unbounded $*$-representation. 
One can easily define a ``well-behaved'' class of such unbounded 
representations and classify them up to unitary equivalence using the same 
technique.

In the paper we use the following standard notations:  ${\mathbb R}$ is the 
set of real numbers, ${\mathbb R}^+$ is the set of nonnegative real numbers,
${\mathbb Z}$ denotes the set of integers, 
${\mathbb Z}^+=\{0,1,2,\ldots\}$. 
\bigskip

\noindent
{\bf 1. The $*$-algebra $Pol(Mat_{2,2})_q$ and its $*$-representations}
\bigskip

Let $q\in (0,1)$. The $*$-algebra
$Pol(Mat_{2,2})_q$, a  $q$-analogues of polynomials on the space 
 $Mat_{2,2}$  of complex $2\times 2$ matrices, is given by its generators  
$\{z_{a}^{\alpha}\}_{a=1,2;\alpha=1,2}$ and the following commutation
relations:
\begin{equation}
\begin{array}{rclrcl}\label{ff1}
z_1^1z_2^1&=&qz_2^1z_1^1,& z_2^1z_1^2&=&z_1^2z_2^1,\\
z_1^1z_1^2&=&qz_1^2z_1^1,& z_2^1z_2^2&=&qz_2^2z_2^1,\\
z_1^1z_2^2-z_2^2z_1^1&=&(q-q^{-1})z_1^2z_2^1,& z_1^2z_2^2&=&qz_2^2z_1^2,\\
\end{array}
\end{equation}

\begin{equation}
\begin{array}{rcl}\label{ff2}
(z_1^1)^*z_1^1&=&q^2z_1^1(z_1^1)^*-
(1-q^2)(z_2^1(z_2^1)^*+z_1^2(z_1^2)^*)+\\
&&\qquad\quad+q^{-2}(1-q^2)^2z_2^2(z_2^2)^*+1-q^2,\\
(z_2^1)^*z_2^1&=&q^2z_2^1(z_2^1)^*-(1-q^2)z_2^2(z_2^2)^*+1-q^2,\\
(z_1^2)^*z_1^2&=&q^2z_1^2(z_1^2)^*-(1-q^2)z_2^2(z_2^2)^*+1-q^2,\\
(z_2^2)^*z_2^2&=& q^2z_2^2(z_2^2)^*+1-q^2,
\end{array}
\end{equation}

\begin{eqnarray}
&\mbox{and}\nonumber\\
&\begin{array}{rclrcl}\label{ff3}
(z_1^1)^*z_2^1-qz_2^1(z_1^1)^*&=&(q-q^{-1})z_2^2(z_1^2)^*,& 
(z_2^2)^*z_2^1&=&qz_2^1(z_2^2)^*,\\
(z_1^1)^*z_1^2-qz_1^2(z_1^1)^*&=&(q-q^{-1})z_2^2(z_2^1)^*,& 
(z_2^2)^*z_1^2&=&qz_1^2(z_2^2)^*,\\
(z_1^1)^*z_2^2&=&z_2^2(z_1^1)^*,& 
(z_2^1)^*z_1^2&=&z_1^2(z_2^1)^*.
\end{array}
\end{eqnarray}

Consider a representation $\pi$ of $Pol(Mat_{2,2})_q$ on a separable 
Hilbert space $H$ by bounded operators.  
The theorem below gives the complete classification of such irreducible 
representations up to a unitary equivalence.
\begin{theorem}
Any irreducible representation $\pi$ is unitarily equivalent to one from the 
following 6 series:

1) one-dimensional representations $\xi_{\varphi_1,\varphi_2}$ 

\begin{equation}\label{f1}
\xi_{\varphi_1,\varphi_2}(z_1^1)=q^{-1}e^{i\varphi_1},\quad 
\xi_{\varphi_1,\varphi_2}(z_2^1)=\xi_{\varphi_1,\varphi_2}(z_1^2)=0,\quad 
\xi_{\varphi_1,\varphi_2}(z_2^2)=e^{i\varphi_2},
\end{equation}
$\varphi_i\in[0,2\pi)$;

2) infinite-dimensional representations $\pi_{\varphi}$ on $H=l_2({\mathbb Z}^+)$
\begin{eqnarray}\label{f2}
&&\begin{array}{rcl}
\pi_{\varphi}(z_1^1)e_{k}&=&q^{-1}\sqrt{1-q^{2(k+1)}}e_{k+1},\\
\pi_{\varphi}(z_2^2)e_{k}&=&e^{i\varphi}e_{k},
\end{array}\\
&&\pi_{\varphi}(z_2^1)=\pi_{\varphi}(z_1^2)=0,\nonumber
\end{eqnarray}
$\varphi\in[0,2\pi)$;

3)  infinite-dimensional representations $\rho_{\varphi_1,\varphi_2}$ on
$H=l_2({\mathbb Z}^+)$
\begin{eqnarray}\label{f3}
\begin{array}{rcl}
\rho_{\varphi_1,\varphi_2}(z_1^1)e_{k}&=&-e^{i(\varphi_1+\varphi_2)}q^{-1}\sqrt{1-q^{2k}}e_{k-1},\\
\rho_{\varphi_1,\varphi_2}(z_2^1)e_{k}&=&e^{i\varphi_1}q^ke_{k},\\
\rho_{\varphi_1,\varphi_2}(z_1^2)e_{k}&=&e^{i\varphi_2}q^ke_{k},\\
\rho_{\varphi_1,\varphi_2}(z_2^2)e_{k}&=&\sqrt{1-q^{2(k+1)}}e_{k+1},
\end{array}
\end{eqnarray}
$\varphi_i\in[0,2\pi)$;

4a)  infinite-dimensional representations $\rho_{\varphi}^1$ on 
$H=l_2({\mathbb Z}^+\times{\mathbb Z}^+)$

\begin{eqnarray}\label{f4a}
\begin{array}{rcl}
\rho_{\varphi}^1(z_1^1)e_{m,k}&=&-e^{i\varphi}q^{-1}\sqrt{1-q^{2(m+1)}}\sqrt{1-q^{2k}}e_{m+1,k-1},\\
\rho_{\varphi}^1(z_2^1)e_{m,k}&=&q^k\sqrt{1-q^{2(m+1)}}e_{m+1,k},\\
\rho_{\varphi}^1(z_1^2)e_{m,k}&=&e^{i\varphi}q^ke_{m,k},\\
\rho_{\varphi}^1(z_2^2)e_{m,k}&=&\sqrt{1-q^{2(k+1)}}e_{m,k+1},
\end{array}
\end{eqnarray}
$\varphi\in[0,2\pi)$;

4b) infinite-dimensional representations $\rho_{\varphi}^2$ on
$H=l_2({\mathbb Z}^+\times{\mathbb Z}^+)$

\begin{eqnarray}\label{f4b}
\begin{array}{rcl}
\rho_{\varphi}^2(z_1^1)e_{m,k}&=&-e^{i\varphi}q^{-1}\sqrt{1-q^{2(m+1)}}\sqrt{1-q^{2k}}e_{m+1,k-1},\\
\rho_{\varphi}^2(z_2^1)e_{m,k}&=&e^{i\varphi}q^ke_{m,k},\\
\rho_{\varphi}^2(z_1^2)e_{m,k}&=&q^k\sqrt{1-q^{2(m+1)}}e_{m+1,k},\\
\rho_{\varphi}^2(z_2^2)e_{m,k}&=&\sqrt{1-q^{2(k+1)}}e_{m,k+1},
\end{array}
\end{eqnarray}
$\varphi\in[0,2\pi)$;

5) infinite-dimensional representations $\hat\rho_{\varphi}$ on 
$H=l_2({\mathbb Z}^+\times{\mathbb Z}^+
\times{\mathbb Z}^+)$

\begin{eqnarray}\label{f5}
\begin{array}{rcl}
\rho (z_1^1)e_{m,l,k}&=&e^{i\varphi}q^{m+l}e_{m,l,k}-\\
&&-q^{-1}\sqrt{(1-q^{2(l+1)})(1-q^{2(m+1)})(1-q^{2k})}e_{m+1,l+1, k-1},\\
\rho (z_2^1)e_{m,l,k}&=&q^k\sqrt{1-q^{2(m+1)}}e_{m+1,l,k},\\
\rho (z_1^2)e_{m,l,k}&=&q^k\sqrt{1-q^{2(l+1)}}e_{m,l+1,k},\\
\rho (z_2^2)e_{m,l,k}&=&\sqrt{1-q^{2(k+1)}}e_{m,l,k+1},
\end{array}
\end{eqnarray}
$\varphi\in[0,2\pi)$;

6) infinite-dimensional representation $\rho$ on $H=l_2({\mathbb Z}^+\times{\mathbb Z}^+
\times{\mathbb Z}^+\times{\mathbb Z}^+)$

\begin{eqnarray}\label{f6}
\begin{array}{rcl}
\rho (z_1^1)e_{s,m,l,k}&=&q^{m+l}\sqrt{1-q^{2(s+1)}}e_{s+1,m,l,k}-\\
&&-q^{-1}\sqrt{(1-q^{2(l+1))}(1-q^{2(m+1)})(1-q^{2k})}e_{s,m+1,l+1, k-1},\\
\rho (z_2^1)e_{s,m,l,k}&=&q^k\sqrt{1-q^{2(m+1)}}e_{s,m+1,l,k},\\
\rho (z_1^2)e_{s,m,l,k}&=&q^k\sqrt{1-q^{2(l+1)}}e_{s,m,l+1,k},\\
\rho (z_2^2)e_{s,m,l,k}&=&\sqrt{1-q^{2(k+1)}}e_{s,m,l,k+1}.
\end{array}
\end{eqnarray}

\end{theorem}
 
\begin{proof}
Let us consider a $*$-subalgebra $\mathcal B$ of $Pol(Mat_{2,2})_q$ which is
generated by 
$z_2^1$, $z_1^2$, $z_2^2$ and  
 $(z_2^1)^*$, $(z_1^2)^*$, $(z_2^2)^*$.
Direct computation shows that $z_2^1(z_2^1)^*$, $z_1^2(z_1^2)^*$, 
$z_2^2(z_2^2)^*$ generate a commutative $*$-subalgebra of $\cal B$ and 
satisfy the 
following relations:
\begin{equation}\label{r}
(z_a^{\alpha}(z_a^{\alpha})^*)z_b^{\beta}=z_b^{\beta}F_{ba}^{\beta\alpha}
(z_2^1(z_2^1)^*, z_1^2(z_1^2)^*, z_2^2(z_2^2)^*)
\end{equation}
where 
\begin{multline*}
{\mathbb F}_{21}(x_1,x_2,x_3)=(F_{22}^{11}(x_1,x_2,x_3), 
F_{21}^{12}(x_1,x_2,x_3), F_{22}^{12}(x_1,x_2,x_3))=\\
=(q^2x_1-(1-q^2)(x_3-1), x_2,x_3),\\
{\mathbb F}_{12}(x_1,x_2,x_3)=(F_{12}^{21}(x_1,x_2,x_3), 
F_{11}^{22}(x_1,x_2,x_3), 
F_{12}^{22}(x_1,x_2,x_3))=\\
=(x_1,q^2x_2-(1-q^2)(x_3-1),x_3),\\
{\mathbb F}_{22}(x_1,x_2,x_3)=(F_{22}^{21}(x_1,x_2,x_3), 
F_{21}^{22}(x_1,x_2,x_3), F_{22}^{22}(x_1,x_2,x_3))=\\
=(q^2x_1, q^2x_2, q^2(x_3-1)+1).
\end{multline*}
The functions ${\mathbb F}_{21}$, ${\mathbb F}_{12}$, 
${\mathbb F}_{22}:{\mathbb R}^3\rightarrow {\mathbb R}^3$ define an action
of ${\mathbb Z}^3$ on ${\mathbb R}^3$ with orbits
\begin{multline*}
\Omega_{x_1,x_2,x_3}=\{{\mathbb F}_{21}^{(m)}({\mathbb F}_{12}^{(l)}
({\mathbb F}_{22}^{(k)}(x_1,x_2,x_3)))=\\
\qquad =(q^{2k}(q^{2m}x_1-(1-q^{2m})(x_3-1)), 
q^{2k}(q^{2l}x_2-(1-q^{2l})(x_3-1)), \\q^{2k}(x_3-1)+1), m,l,k\in{\mathbb Z}\}.
\end{multline*}

Here and in the sequel we denote by ${\mathbb F}_{a\alpha}^{(m)}$ the $m$-th
iteration of ${\mathbb F}_{a\alpha}$ and $({\mathbb F}_{a\alpha}^{(m)})_i$, 
$i=1,2,3$, the $i$-th coordinate of ${\mathbb F}_{a\alpha}^{(m)}$.
Let $\pi$ be a $*$-representation of $Pol(Mat_{2,2})_q$ on a Hilbert
 space $H$ by bounded operators, let
${\mathbb E}(\cdot)$ be the resolution of the identity for the commutative 
family  ${\mathbb A}_{\pi}$ of the positive 
operators $\pi(z_2^1)\pi(z_2^1)^*$, $\pi(z_1^2)\pi(z_1^2)^*$, 
$\pi(z_2^2)\pi(z_2^2)^*$ 
and let $\sigma_{\pi}$ be the joint spectrum of the family 
${\mathbb A}_{\pi}$.

Next step is to show that any irreducible representation is concentrated on an orbit of this dynamical system.
\begin{lemma}
If $\pi$ is an irreducible representation of $Pol(Mat_{2,2})_q)$ then the 
spectral measure ${\mathbb E}(\cdot)$ is ergodic with respect to  the action 
of the dynamical system generated by 
${\mathbb F}_{21}$, ${\mathbb F}_{12}$, ${\mathbb F}_{22}$ and there 
exists an orbit $\Omega_{x_1,x_2,x_3}$  such that 
${\mathbb E}(\Omega_{x_1,x_2,x_3})=I$.
\end{lemma}
\begin{proof}
From (\ref{r}) and the spectral theorem  it follows that
\begin{eqnarray*}
\begin{array}{rcl}
{\mathbb E}(\Delta)\pi(z_b^{\beta})&=&\pi(z_b^{\beta}){\mathbb E}
({\mathbb F}_{b\beta}^{(-1)}(\Delta)),\\
{\mathbb E}(\Delta)\pi(z_b^{\beta})^*&=&\pi(z_b^{\beta})^*{\mathbb E}
({\mathbb F}_{b\beta}(\Delta)),
\end{array}
\end{eqnarray*}
for any $\Delta\in{\mathfrak B}({\mathbb R}^3)$. Hence any subset $\Delta$ 
such that ${\mathbb F}_{b\beta}^{(-1)}(\Delta)\subseteq\Delta$, 
${\mathbb F}_{b\beta}(\Delta)\subseteq\Delta$, $(b,\beta)=(2,1)$, $(1,2)$, 
$(2,2)$ defines a subspace
${\mathbb E}(\Delta)H$ which is invariant with respect to the operators
$\pi(z_b^{\beta})$, $\pi(z_b^{\beta})^*$ for any $(b,\beta)$ as above.
Moreover,  such subspace is invariant with respect to any 
operator of the representation $\pi$. In fact, the following relations hold in
$Pol(Mat_{2,2})_q$
\begin{eqnarray}\label{equ}
&z_a^{\alpha}(z_a^{\alpha})^*z_1^1=z_1^1z_a^{\alpha}(z_a^{\alpha})^*-
(-1)^{a+\alpha}(q-q^{-1})z_2^1z_1^2(z_2^2)^*
\end{eqnarray}
$(a,\alpha)=(2,1)$, $(1,2)$, 
$(2,2)$, which gives
\begin{eqnarray*}
{\mathbb E}({\mathbb R}^3\setminus\Delta)\pi(z_a^{\alpha}(z_a^{\alpha})^*)
\pi(z_1^1){\mathbb E}(\Delta)={\mathbb E}({\mathbb R}^3\setminus\Delta)
\pi(z_1^1)\pi(z_a^{\alpha}(z_a^{\alpha})^*){\mathbb E}(\Delta)-\\
-(-1)^{a+\alpha}
(q-q^{-1}){\mathbb E}({\mathbb R}^3\setminus\Delta)\pi(z_2^1)\pi(z_1^2)
\pi(z_2^2)^*{\mathbb E}(\Delta)
\end{eqnarray*}
Therefore if $\Delta\in {\mathfrak B}({\mathbb R}^3)$ is invariant with 
respect 
to all ${\mathbb F}_{b\beta}^{(-1)}$ and ${\mathbb F}_{b\beta}$ we obtain
$$\pi(z_a^{\alpha}(z_a^{\alpha})^*){\mathbb E}({\mathbb R}^3\setminus\Delta)
\pi(z_1^1){\mathbb E}(\Delta)=
{\mathbb E}({\mathbb R}^3\setminus\Delta)\pi(z_1^1)
{\mathbb E}(\Delta)\pi(z_a^{\alpha}(z_a^{\alpha})^*)$$
and hence
$${\mathbb E}(\Delta'){\mathbb E}({\mathbb R}^3\setminus\Delta)
\pi(z_1^1){\mathbb E}(\Delta)={\mathbb E}({\mathbb R}^3\setminus\Delta)
\pi(z_1^1){\mathbb E}(\Delta){\mathbb E}(\Delta')$$
for any $\Delta'\in{\mathfrak B}({\mathbb R}^3)$.
Taking $\Delta'=\Delta$ gives ${\mathbb E}({\mathbb R}^3\setminus\Delta)
\pi(z_1^1){\mathbb E}(\Delta)=0$, i.e.  
$\pi(z_1^1){\mathbb E}(\Delta)H\subseteq
{\mathbb E}(\Delta)H$. Similarly, $\pi(z_1^1)^*{\mathbb E}(\Delta)H\subseteq
{\mathbb E}(\Delta)H$.
The ergodicity of the measure ${\mathbb E}(\cdot)$ follows immediately, i.e., 
${\mathbb E}(\Delta)=I$ or $0$ for any Borel $\Delta$ which is invariant
with respect to ${\mathbb F}_{b\beta}$, ${\mathbb F}_{b\beta}^{(-1)}$.

The simplest invariant sets are the  orbits of the dynamical system. 
The next step 
is to show that only atomic measures concentrated on an orbit give rise to 
irreducible representation of the $*$-algebra.
It is easily seen that the dynamical system generated by 
${\mathbb F}_{b\beta}$ is one-to-one and 
possesses a measurable section, i.e., a  set 
$\tau\in{\mathfrak B}({\mathbb R}^3)$ 
which intersects any orbit in a single point. This implies that any
ergodic measure is concentrated on a single orbit of the dynamical system and 
therefore ${\mathbb E}(\Omega_{x_1,x_2,x_3})=I$ for some orbit 
$\Omega_{x_1,x_2,x_3}$.
\end{proof}

We now clarify  which orbits $\Omega_{x_1,x_2,x_3}$ give rise to  bounded
irreducible representation $\pi$, i.e.,  
$\sigma_{\pi}\subseteq\Omega_{x_1,x_2,x_3}$, and
classify all such representations up to unitary equivalence. 

We claim first that there is no bounded representations $\pi$ with   
$\sigma_{\pi}\subseteq\Omega_{x_1,x_2,x_3}$ if $x_3>1$.  From (\ref{r}) we have
\begin{equation}\label{s}
\pi(z_b^{\beta})H_x\subseteq H_{{\mathbb F}_{b\beta}(x)},\quad
\pi(z_b^{\beta})^*H_x\subseteq H_{{\mathbb F}_{b\beta}^{(-1)}(x)},
\end{equation}
 where $H_x$ is the
eigenspace for ${\mathbb A}_{\pi}$ corresponding to the eigenvalue
$x\in{\mathbb R}^3$. Since $y=(y_1,y_2,y_3)\in\Omega_{x_1,x_2,x_3}$, where
$x_3>1$, implies $y_3>1$ we conclude that $\pi(z_2^2)\pi(z_2^2)^*\geq 1$ and 
$\ker\pi(z_2^2)=\ker\pi(z_2^2)^*=\{0\}$. This
clearly forces ${\mathbb F}_{22}^{(k)}(y)\in\sigma_{\pi}$ for any 
$k\in{\mathbb Z}$. However, the set 
$\{{\mathbb F}_{22}^{(k)}(y), k\in{\mathbb Z}\}$ is unbounded which 
contradicts the boundness of the representation $\pi$.
Similar arguments show that there is no bounded representation $\pi$ with
$\sigma_{\pi}\subseteq\Omega_{x_1,x_2,1}$, $x_1\ne 0$ or $x_2\ne 0$. In this case
$\Omega_{x_1,x_2,1}=\{(q^{2(k+m)}x_1,q^{2(k+l)}x_2,1), k,l,m\in {\mathbb Z}\}$.
The only possibility is $\sigma_{\pi}=\Omega_{0,0,1}=\{(0,0,1)\}$ and in this case we obtain 
$\pi(z_2^1)=\pi(z_1^2)=0$, $\pi(z_2^2)\pi(z_2^2)^*=I$. It follows now from 
(\ref{ff1})--(\ref{ff3}) that  $\pi(z_2^2)$, $\pi(z_1^1)$ satisfy the relations
\begin{eqnarray}\label{x}
&\pi(z_1^1)^*\pi(z_1^1)=q^2\pi(z_1^1)\pi(z_1^1)^*+(q^{-2}-1),\nonumber\\
&{}[\pi(z_1^1),\pi(z_1^2)]=0,\quad [\pi(z_1^1)^*,\pi(z_1^2)]=0,\\
&\pi(z_2^2)^*\pi(z_2^2)=\pi(z_2^2)\pi(z_2^2)^*=I.\nonumber
\end{eqnarray}
This implies that $\pi(z_2^2)$ commutes with all images of the generators in
the algebra under th e representation $\pi$ and therefore $\pi(z_2^2)$ is a 
multiple of the identity operator if $\pi$ is irreducible. By (\ref{x}) we have
also $\pi(z_2^2)=e^{i\varphi_2}I$, $\varphi_2\in[0,2\pi)$.
 Irreducible representations of
the relation
$(z_1^1)^*z_1^1=q^2z_1^1(z_1^1)^*+(q^{-2}-1)$ are well-known and can be easily 
calculated using the method of dynamical systems (see \cite[Chapter~2]{OS}). 
Any such representation is  either one-dimensional: 
$\xi_{\varphi_1}(z_1^1)=q^{-1}e^{i\varphi_1}$, $\varphi_1\in [0,2\pi)$, or
infinite-dimensional which is unitary equivalent to the following one
$\pi_{\varphi}(z_1^1)e_k=q^{-1}\sqrt{1-q^{2(k+1)}}e_{k+1}$. The corresponding
irreducible representations of $Pol(Mat_{2,2})_q$ are 
$\xi_{\varphi_1,\varphi_2}$ and $\pi_{\varphi}$.

Since $\sigma_{\pi}\subseteq ({\mathbb R}^+)^3$ and 
$({\mathbb F}_{22}^{(k)})_3(x_1,x_2,x_3)=q^{2k}(x_3-1)+1\rightarrow -\infty$ as
$k\rightarrow-\infty$, it follows from (\ref{s}) that 
$\ker\pi(z_2^2)^*\ne\{0\}$, $\ker\pi(z_2^2)\pi(z_2^2)^*\ne\{0\}$ and the 
corresponding orbit contains a point $(x_1,x_2,0)$. We have
$\Omega_{x_1,x_2,0}=\{(q^{2k}(q^{2m}(x_1-1)+1), q^{2k}(q^{2l}(x_2-1)+1), 
1-q^{2k}), m,l,k\in{\mathbb Z}\}$.
Similar arguments show that $\sigma_{\pi}\subseteq\Omega_{x_1,x_2,0}$, where
$x_1>1$ or $x_2>1$, is impossible if the representation $\pi$ is bounded. From
the positiveness of $\sigma_{\pi}$ we obtain also that the only orbits 
corresponding to irreducible representation of the $*$-algebra are
$\Omega_{1,1,0}$, $\Omega_{1,0,0}$, $\Omega_{0,1,0}$, $\Omega_{0,0,0}$ and
$\Omega_{0,0,1}$ The last one was treated above.  

We consider now the case $\sigma_{\pi}\subseteq\Omega_{x_1,x_2,x_3}$, $x_3=0$.
Let $P_y$, $y=(y_1,y_2,y_3)$ be the projection onto the eigenspace 
corresponding to
the eigenvalue $y$. Using  (\ref{equ}) we get
$$(z_k-y_k)P_z\pi(z_1^1)P_y=\pm (q-q^{-1})P_z\pi(z_2^1)\pi(z_1^2)
\pi(z_2^2)^*P_y$$
(``+'' for $k=1,2$ and ``-'' for $k=3$)  $z,y\in{\mathbb R}^3$.
By (\ref{s}) we have $\pi(z_2^1)\pi(z_1^2)\pi(z_2^2)^*H_y\subseteq
H_{{\mathbb F}_{21}({\mathbb F}_{12}({\mathbb F}_{22}^{(-1)}(y)))}$ and
$$\pi(z_2^1)\pi(z_1^2)\pi(z_2^2)^*P_y=P_{{\mathbb F}_{21}({\mathbb F}_{12}
({\mathbb F}_{22}^{(-1)}(y)))}\pi(z_2^1)\pi(z_1^2)\pi(z_2^2)^*P_y.$$
Setting  $P_{m,l,k}$  the projection onto an eigenspace which 
corresponds to the eigenvalue
${\mathbb F}_{21}^{(m)}({\mathbb F}_{12}^{(l)}({\mathbb F}_{22}^{(m)}
(x_1,x_2,0)))$ we obtain 
$$\pi(z_1^1)P_{m,l,k}=P_{m,l,k}\pi(z_1^1)P_{m,l,k}+ 
P_{m+1,l+1,k-1}\pi(z_1^1)P_{m,l,k},$$ i.e., 
$$\pi(z_1^1)H_{m,l,k}\subseteq H_{m,l,k}\oplus H_{m+1,l+1,k-1}.$$
Moreover, $P_{m+1,l+1,k-1}\pi(z_1^1)P_{m,l,k}=-q^{1-2k}\pi(z_2^1)\pi(z_1^2)\pi(z_2^2)^*P_{m,l,k}$. The operator $\pi(z_1^1)$ can be written  as a sum of
its diagonal part $\pi(z_1^1)_0=\sum_{m,l,k}P_{m,l,k}\pi(z_1^1)P_{m,l,k}$, and 
the operator
$-\sum_{m,l,k}q^{1-2k}\pi(z_2^1)\pi(z_1^2)\pi(z_2^2)^*P_{m,l,k}=
-q\pi(z_2^1)\pi(z_1^2)\pi(z_2^2)^*(1-\pi(z_2^2(z_2^2)^*))^{-1}$.

Let now $\sigma_{\pi}\subseteq\Omega_{x_1,x_2,0}$, where
$x_1\ne 0$ or $x_2\ne 0$.
It follows from (\ref{ff1})-(\ref{ff3}) by direct computation that 
$$\pi(z_1^1)_0^*\pi(z_1^1)_0=
q^2\pi(z_1^1)_0\pi(z_1^1)_0^*.$$ The only bounded operator $\pi(z_1^1)_0$ 
satisfying this relation is the zero-operator. Therefore
$$\pi(z_1^1)=
-q\pi(z_2^1)\pi(z_1^2)\pi(z_2^2)^*(1-\pi(z_2^2(z_2^2)^*))^{-1}$$
and $\pi$ is irreducible iff so is the family
$(\pi(z_2^1), \pi(z_1^2), \pi(z_2^2), \pi(z_2^1)^*, \pi(z_1^2)^*, \pi(z_2^2)^*)$.
Let $\pi(z_a^{\alpha})=U_a^{\alpha}\sqrt{\pi(z_a^{\alpha})^*\pi(z_a^{\alpha})}$
be the polar decomposition of $\pi(z_a^{\alpha})$.
Using easy arguments  one can show that $[U_a^{\alpha},U_b^{\beta}]=
[(U_a^{\alpha})^*,U_b^{\beta}]=0$, $(a,\alpha)\ne(b,\beta)$ and   
$$(z_a^{\alpha}(z_a^{\alpha})^*)(U_b^{\beta})=(U_b^{\beta})F_{ba}^{\beta\alpha}
(z_2^1(z_2^1)^*, z_1^2(z_1^2)^*, z_2^2(z_2^2)^*).$$
Here $(a,\alpha)$, $(b,\beta)\in\{(1,2), (2,1),(2,2)\}$.
Moreover, if $\sigma_{\pi}\subseteq\Omega_{1,x_2,0}$ 
($\sigma_{\pi}\subseteq\Omega_{x_1,1,0}$) we have $U_2^1$ ($U_1^2$ 
respectively)  commutes
with any operators from the family ${\mathbb A}_{\pi}$ and therefore with any 
operator of the representation. This clearly forces $U_2^1=e^{i\varphi_1}I$,
$\varphi_1\in [0,2\pi)$
($U_1^2=e^{i\varphi_2}I$, $\varphi_2\in[0,2\pi)$ respectively). 
Let $\sigma_{\pi}\subseteq \Omega_{0,1,0}$. Consider
 $e_{k,l}=(U_2^1)^k(U_2^2)^le$, 
$e\in\ker\pi(z_2)\pi(z_2^2)^*\cap\ker\pi(z_2^1)\pi(z_2^1)^*$, 
$k,l\in{\mathbb Z}^+$. Then $\{e_{k,l}, k,l\in{\mathbb Z}^+\}$ is an 
orthonormal system which defines an invariant subspace.
The corresponding irreducible representation is $\rho_{\varphi}^2$.
 Analogously $(U_1^2)^k(U_2^2)^le=e_{k,l}$, 
$e\in\ker\pi(z_2)\pi(z_2^2)^*\cap\ker\pi(z_1^2)\pi(z_1^2)^*$, 
$k,l\in{\mathbb Z}^+$,
build an orthonormal basis of an irreducible representation space if
$\sigma_{\pi}\subseteq\Omega_{1,0,0}$, the corresponding action is given by 
formulae (\ref{f4b}).
If $\sigma_{\pi}\subseteq\Omega_{1,1,0}$ we have that  
$l.s.\{(U_2^2)^ke=e_k, k\in{\mathbb Z}^+\}$, 
$e\in\ker\pi(z_2^2)\pi(z_2^2)^*$ is invariant with the corresponding action 
given by (\ref{f3}).
 
We now turn to the case $\sigma_{\pi}\subset\Omega_{0,0,0}$.
From (\ref{ff1})--(\ref{ff3}) we have
\begin{eqnarray}
\pi(z_1^1)_0\pi(z_1^2)=q\pi(z_1^2)\pi(z_1^1)_0\quad \pi(z_1^1)_0^*\pi(z_1^2)=
q\pi(z_1^2)\pi(z_1^1)_0^*,\nonumber\\
\pi(z_1^1)_0\pi(z_2^1)=q\pi(z_2^1)\pi(z_1^1)_0\quad \pi(z_1^1)_0^*\pi(z_2^1)=
q\pi(z_2^1)\pi(z_1^1)_0^*,\nonumber\\
\pi(z_1^1)_0\pi(z_2^2)=\pi(z_2^2)\pi(z_1^1)_0\quad \pi(z_1^1)_0^*\pi(z_2^2)=
\pi(z_2^2)\pi(z_1^1)_0^*,\nonumber
\end{eqnarray}
$$\pi(z_1^1)_0^*\pi(z_1^1)_0P_{m,l,k}=
q^2\pi(z_1^1)_0\pi(z_1^1)_0^*P_{m,l,k}+(1-q^2)q^{2(m+l)}P_{m,l,k}.$$
Note that $\pi(z_1^1)_0P_{m,l,k}H\subseteq P_{m,l,k}H$, 
$\pi(z_1^1)_0^*P_{m,l,k}H\subseteq P_{m,l,k}H$.
Moreover, it follows from the above relation  that
if $\pi$ is irreducible then the family
$(\pi(z_1^1)_0,\pi(z_1^1)^*_0)$ restricted to the subspace $P_{m,l,k}H$ is
irreducible for any $m,l,k\in{\mathbb Z}^+$.
We have  
$$a^*a=q^2aa^*+(1-q^2),$$
where $a=\pi(z_1^1)_0P_{0,0,0}$
Any  irreducible family ($a$, $a^*$) is either one-dimensional and given by
$a=e^{i\varphi}$, $\varphi\in [0,2\pi)$, or infinite dimensional defined  on 
$l_2({\mathbb Z}^+)$
 by 
$ae_{s}=\sqrt{1-q^{2(s+1)}}e_{s+1}$.
These representations give rise to irreducible representations of the 
$*$-algebra $Pol(Mat_{2,2})_q$. Namely, in the first case we have
that  $e_{m,l,k}=(U_2^1)^m(U_1^2)^l(U_2^2)^ke$, where $e\in P_{0,0,0}H=
\ker\pi(z_2^2)\pi(z_2^2)^*\cap \ker\pi(z_1^2)\pi(z_1^2)^*\cap 
\ker\pi(z_2^1)\pi(z_2^1)^*$, $m,l,k\in{\mathbb Z}^+$, define an orthonormal 
basis of the space where the irreducible representation $\hat\rho_{\varphi}$ 
acts, and for the second irreducible family we have that 
 $e_{s,m,l,k}=(U_2^1)^m(U_1^2)^l(U_2^2)^ke_s$, $s,m,l,k\in{\mathbb Z}^+$,
define an orthonormal basis of the space where the irreducible representation
$\rho$ acts. This finishes the proof.
\end{proof}

{\it Comments.} It follows from the proof that for any representation $\pi$ on
a Hilbert space $H_{\pi}$
the family  of self-adjoint operators $\pi(z_2^2(z_2^2)^*)$, 
$\pi(z_2^1(z_2^1)^*)$, $\pi(z_1^2(z_1^2)^*)$, 
$\pi(z_1^1)_0\pi(z_1^1)^*_0$,  where
$$ \pi(z_1^1)_0=\pi(z_1^1)-\left\{\begin{array}{ll}
0,&  \pi(z_2^2(z_2^2)^*)=I\\
-q\pi(z_2^1)\pi(z_1^2)\pi(z_2^2)^*(1-\pi(z_2^2(z_2^2)^*)^{-1},& 
\pi(z_2^2(z_2^2)^*)\ne I\end{array}\right.$$ generates a commutative 
$*$-subalgebra $\mathbb A$ in $B(H_{\pi})$, the bounded operators on $H_{\pi}$.
Moreover, any irreducible representation of $Pol(Mat_{2,2})_q$ is a weight 
representation with respect to this algebra, i.e., ${\mathbb A}$ can be 
diagonalized, and the spectrum of  ${\mathbb A}$ is simple.

A  question which arise here is how to generalise the method to higher 
dimension matrix balls and classify $*$-representations of the corresponding 
$*$-algebras. In principle, just analysing the commutation relations between 
the generators in the $*$-algebra one can  find a commutative $*$-subalgebra
of $Pol(Mat_{m,n})_q$ or some its localisation and show that  any irreducible 
representation $\pi$ is a weight representation
with respect to this commutative $*$-algebra having a  simple spectrum in 
this representation. However, the computations can be extremely difficult in 
general.

\begin{remark}\rm
The polynomial algebra on the vector space $Mat_{2,2}$ can be supplied with a 
Poisson structure. Writing $q=e^{-h}$ we have that  
$Pol(Mat_{2,2})_{exp(-h)}$ is
an associative algebra over the ring of formal series ${\mathbb C}[[h]]$  and
$$Pol(Mat_{2,2})\simeq Pol(Mat_{2,2})_{exp(-h)}/hPol(Mat_{2,2})_{exp(-h)}.$$ 
The Poisson bracket now is given by
$$\{a\mod h,b\mod h\}=-ih^{-1}(ab-ba)\mod h$$ 
for any $a$, $b\in Pol(Mat_{2,2})_{exp(-h)}$.
The problem now is to define the symplectic leaves of this Poisson structure.
Any primitive ideal $\ker\pi$, where $\pi$ is an irreducible representation
of $Pol(Mat_{2,2})_q$,  defines a maximal Poisson ideal $I_{\pi}=\ker\pi\mod h$
 of the algebra
$Pol(Mat_{2,2})$ ordered by inclusion and hence the closure of a 
symplectic leaf which is given by 
$\{x\in Mat_{2,2}\mid f(x)=0,f\in I_{\pi}\}$.
As in the case of ${\mathbb C}(SU(n))_q$ (see \cite{SoV}) one can expect 
that there is a one-to-one correspondence between irreducible representations
 (bounded irreducible representations) of 
$Pol(Mat_{2,2})_q$ and symplectic leaves (bounded symplectic leaves) in 
$Mat_{2,2}$.
\end{remark}

\vspace{0.1cm}
{\bf Acknowledgement.}
I am sincerely grateful to Prof. L. Vaksman for having communicated the 
problem to me and for a number helpful discussions during the preparation of 
the paper.

\end{document}